\documentclass[12pt,reqno]{amsart}
\usepackage{amsfonts}
\usepackage{amsmath,amsthm,amsfonts,amssymb}
\newtheorem{lemma}{Lemma}[section]
\newtheorem{theorem}{Theorem}[section]
\newtheorem{corollary}{Corollary}[section]

\def\bl{\begin{lemma}}
\def\bt{\begin{theorem}}
\def\el{\end{lemma}}
\def\et{\end{theorem}}
\def\bp{\begin{proof}}
\def\ep{\end{proof}}
\def\bc{\begin{corollary}}
\def\ec{\end{corollary}}

\def\mb{\mathbb}

\def\-{\setminus}

\def\vp{\varphi}

\def\lt{\left}
\def\rt{\right}
\def\+{\bigcup}
\def\.{\bigcap}

\def\ll{\langle}
\def\rl{\rangle}

\title[]
{A right inverse of differential operator $\Delta+a$ in weighted Hilbert space $L^2(\mb{R}^n,e^{-|x|^2})$}

\author []{Shaoyu Dai$^1$, Yang Liu$^2$ and Yifei Pan$^3$}

\address{1 School of Mathematics, Southeast University, Nanjing, 210096, China;}
\address{
Department of Mathematics, Jinling Institute of Technology, Nanjing, 211169, China.}
\address{\it E-mail address: dymdsy@163.com}

\address{2 Department of Mathematics, Zhejiang Normal University, Jinhua, 321004, China.}
\address{\it E-mail address: liuyang@zjnu.edu.cn}

\address{3 Department of Mathematical Sciences, Purdue University Fort Wayne, Fort Wayne, 46805-1499, USA.}
\address{\it E-mail address: pan@pfw.edu}

\begin{document}
\begin{abstract} In this note, we prove
 the existence of weak solutions of a Poisson type equation in the weighted Hilbert space $L^2(\mb{R}^n,e^{-|x|^2})$.
\end{abstract}
\maketitle

\section{Introduction}

In this note, we study the right inverse of differential operator $\Delta+a$ in a Hilbert space by proving the following result on the existence of (global) weak solutions of a Poisson type equation in the weighted Hilbert space $L^2(\mb{R}^n,e^{-|x|^2})$, where $a$ is a real constant throughout.

\bt\label{thk2}
For each $f\in L^2(\mb{R}^n,e^{-|x|^2})$, there exists a weak solution $u\in L^2(\mb{R}^n,e^{-|x|^2})$ solving the equation
$$\Delta u+au=f$$ in $\mb{R}^n$
with the norm estimate $$\int_{\mb{R}^n} u^2e^{-|x|^2}dx\leq
\frac{1}{8n}\int_{\mb{R}^n} f^2e^{-|x|^2}dx.$$
\et

The novelty of Theorem \ref{thk2} is that the differential operator $\Delta+a$ has a bounded right inverse
\begin{align*}
Q: L^2(\mb{R}^n,e^{-|x|^2})&\longrightarrow L^2(\mb{R}^n,e^{-|x|^2}), \\
(\Delta+a)Q&=I
\end{align*}
with the norm estimate $\|Q\|\leq\frac{1}{\sqrt{8n}}$.
In particular, the Laplace operator $\Delta$ has a bounded right inverse $Q_0: L^2(\mb{R}^n,e^{-|x|^2})\longrightarrow L^2(\mb{R}^n,e^{-|x|^2})$, which, to the best of our knowledge, appears to be new.

As a simple consequence of Theorem \ref{thk2}, we can obtain the following
result on the existence of global weak solutions of a Poisson type equation for square integrable functions. This result seems to be classical and well-known, but we could not find an exact reference. At least it is well-known \cite{2} that when $f$ is compact supported in $\mb{R}^n$, then a solution of $-\Delta u=f$ is $\Phi\ast f$, where $\Phi$ is the fundamental solution of Laplace's equation.

\bt\label{th1}
For each $f\in L^2(\mb{R}^n)$ or $f\in L^\infty(\mb{R}^n)$,
there exists a weak solution $u\in L^2_{loc}(\mb{R}^n)$ solving the equation
$$\Delta u+au=f.$$
In particular, the Poisson equation $\Delta u=f$
has a weak solution $u\in L^2_{loc}(\mb{R}^n)$ for $f\in L^2(\mb{R}^n)$ or $f\in L^\infty(\mb{R}^n)$.
\et

The proof of Theorem \ref{th1} follows from the observation that $L^2(\mb{R}^n)\subset L^2(\mb{R}^n,e^{-|x|^2})$, $L^\infty(\mb{R}^n)\subset L^2(\mb{R}^n,e^{-|x|^2})$ and $L^2(\mb{R}^n,e^{-|x|^2})\subset L^2_{loc}(\mb{R}^n)$.

The method employed in this note was motivated from the H\"{o}rmander $L^2$ method \cite{1} for Cauchy-Riemann
equations from several complex variables.

\section{Several lemmas}
Here, we consider weighted Hilbert space
$$L^2(\mb{R}^n,e^{-\vp})
=\{f\mid f\in L^2_{loc}(\mb{R}^n); \int_{\mb{R}^n}f^2 e^{-\vp}dx<+\infty\},$$
where $\vp$ is a nonnegative function on $\mb{R}^n$.
We denote
the weighted inner product for $f,g\in L^2(\mb{R}^n,e^{-\vp})$ by
$\langle f,g\rangle_\vp=\int_{\mb{R}^n}fg e^{-\vp}dx$
and the weighted norm of $f\in L^2(\mb{R}^n,e^{-\vp})$ by
$\|f\|_\vp=\sqrt{\langle f,f\rangle_\vp}.$
Let $k$ be a positive integer, $C^k(\mb{R}^n)$ denote the set of all $k$-times continuously differentiable real-valued functions on $\mb{R}^n$, and $C_0^\infty(\mb{R}^n)$ denote the set of all smooth functions $\phi: \mb{R}^n\rightarrow\mb{R}$ with compact support.
Let $\Delta=\sum^n_{j=1}\frac{\partial^2}{\partial x_j^2}$ be the Laplace operator on $\mb{R}^n$. For $u,f\in L^2_{loc}(\mb{R}^n)$, we say that $f$ is the Laplace of $u$ in the weak sense, written $\Delta u=f$, provided $\int_{\mb{R}^n}u\Delta\phi dx=\int_{\mb{R}^n}f\phi dx$ for all test functions $\phi\in C_0^\infty(\mb{R}^n)$.

Let $\vp\in C^2(\mb{R}^n)$. For $\forall\phi\in C_0^\infty(\mb{R}^n)$, we first define the following formal adjoint of $\Delta$ with respect to the weighted inner product in $L^2\lt(\mb{R}^n,e^{-\vp}\rt)$. Let $u\in L^2_{loc}(\mb{R}^n)$. We calculate as follows.
\begin{align*}
\lt\ll \phi,\Delta u\rt\rl_\vp&=\int_{\mb{R}^n}\phi \Delta u e^{-\vp}dx \\
&=\int_{\mb{R}^n}u \Delta\lt(\phi e^{-\vp}\rt)dx\\
&=\int_{\mb{R}^n}e^{\vp}u \Delta\lt(\phi e^{-\vp}\rt)e^{-\vp}dx\\
&=\lt\ll e^{\vp}\Delta\lt(\phi e^{-\vp}\rt),u\rt\rl_\vp\\
&=:\lt\ll \Delta_\vp^{*}\phi, u\rt\rl_\vp,
\end{align*}
where $\Delta_\vp^{*}\phi =e^{\vp}\Delta\lt(\phi e^{-\vp}\rt)$ is the formal adjoint of $\Delta$ with domain in $C_0^\infty(\mb{R}^n)$. Let $\lt(\Delta+a\rt)_\vp^{*}$ be the formal adjoint of $\Delta+a$ with domain in $C_0^\infty(\mb{R}^n)$. Note that $I_\vp^{*}=I$, where $I$ is  the identity operator. Then $\lt(\Delta+a\rt)_\vp^{*}=\Delta_\vp^{*}+a$.

Let $\nabla$ be the gradient operator on $\mb{R}^n$.
Now we give several key lemmas based on functional analysis.

\bl\label{lemmaifif}
Let $\vp\in C^2(\mb{R}^n)$. For each $f\in L^2(\mb{R}^n,e^{-\vp})$, there exists a global weak solution $u\in L^2(\mb{R}^n,e^{-\vp})$ solving the equation
$$\Delta u+au=f$$ in $\mb{R}^n$
with the norm estimate
$$\|u\|^2_\vp\leq c$$
 if and only if
$$|\langle f,\phi\rangle_\vp|^2\leq c\lt\|\lt(\Delta+a\rt)^*_\vp\phi\rt\|^2_\vp, \ \ \forall\phi\in C_0^\infty(\mb{R}^n),$$
 where $c$ is a constant.
\el
\bp
Let $\Delta+a=H$. Then $\lt(\Delta+a\rt)^*_\vp=H^*_\vp$.

(Necessity) For $\forall\phi\in C_0^\infty(\mb{R}^n)$, from the definition of $H^*_\vp$ and Cauchy-Schwarz inequality, we have
$$|\langle f,\phi\rangle_\vp|^2=|\langle Hu,\phi\rangle_\vp|^2=\lt|\lt\langle u,H^*_\vp\phi\rt\rangle_\vp\rt|^2\leq\|u\|^2_\vp
\lt\|H^*_\vp\phi\rt\|^2_\vp\leq c\lt\|H^*_\vp\phi\rt\|^2_\vp
=c\lt\|\lt(\Delta+a\rt)^*_\vp\phi\rt\|^2_\vp.$$

(Sufficiency) Consider the subspace
$$E=\lt\{H^*_\vp\phi\mid\phi\in C_0^\infty(\mb{R}^n)\rt\}\subset L^2(\mb{R}^n,e^{-\vp}).$$
Define a linear functional $L_f: E\rightarrow\mb{R}$ by
$$L_f\lt(H^*_\vp\phi\rt)=\langle f,\phi\rangle_\vp=\int_{\mb{R}^n}f\phi e^{-\vp}dx.$$
Since
$$\lt|L_f\lt(H^*_\vp\phi\rt)\rt|=\lt|\langle f,\phi\rangle_\vp\rt|
\leq\sqrt{c}\lt\|H^*_\vp\phi\rt\|_\vp,$$
then $L_f$ is a bounded functional on $E$. Let $\overline{E}$ be the closure of $E$ with respect to the norm $\|\cdot\|_\vp$ of $L^2(\mb{R}^n,e^{-\vp})$. Note that $\overline{E}$ is a Hilbert subspace of
$L^2(\mb{R}^n,e^{-\vp})$. So by Hahn-Banach's extension theorem, $L_f$ can be extended to a linear functional $\widetilde{L}_f$ on $\overline{E}$
such that
\begin{equation}\label{29}
\lt|\widetilde{L}_f(g)\rt|\leq\sqrt{c}\lt\|g\rt\|_\vp, \ \ \forall g\in \overline{E}.
\end{equation}
Using the Riesz representation theorem for $\widetilde{L}_f$, there exists a unique $u_0\in \overline{E}$ such that
\begin{equation}\label{30}
\widetilde{L}_f(g)=\langle u_0,g\rangle_\vp, \ \ \forall g\in \overline{E}.
\end{equation}

Now we prove $\Delta u_0+au_0=f$. For $\forall\phi\in C_0^\infty(\mb{R}^n)$, apply $g=H^*_\vp\phi$ in (\ref{30}). Then
$$\widetilde{L}_f\lt(H^*_\vp\phi\rt)=\lt\langle u_0,H^*_\vp\phi\rt\rangle_\vp=\lt\langle Hu_0,\phi\rt\rangle_\vp.$$
Note that
$$\widetilde{L}_f\lt(H^*_\vp\phi\rt)=L_f\lt(H^*_\vp\phi\rt)=\langle f,\phi\rangle_\vp.$$
Therefore,
$$\lt\langle Hu_0,\phi\rt\rangle_\vp=\langle f,\phi\rangle_\vp, \ \ \forall\phi\in C_0^\infty(\mb{R}^n),$$
i.e.,
$$\int_{\mb{R}^n} Hu_0\phi e^{-\vp}dx=\int_{\mb{R}^n} f\phi e^{-\vp}dx, \ \ \forall\phi\in C_0^\infty(\mb{R}^n),$$
i.e.,
$$\int_{\mb{R}^n} (Hu_0-f)\phi e^{-\vp}dx=0, \ \ \forall\phi\in C_0^\infty(\mb{R}^n).$$
Thus, $Hu_0=f$, i.e., $\Delta u_0+au_0=f$.

Next we give a bound for the norm of $u_0$. Let $g=u_0$ in(\ref{29}) and (\ref{30}). Then we have
$$\|u_0\|^2_\vp=\lt|\langle u_0,u_0\rangle_\vp\rt|=\lt|\widetilde{L}_f(u_0)\rt|
\leq\sqrt{c}\lt\|u_0\rt\|_\vp.$$
Therefore, $\|u_0\|¡ª_\vp^2\leq c$.

Note that $u_0\in\overline{E}$ and $\overline{E}\subset L^2(\mb{R}^n,e^{-\vp})$. Then $u_0\in L^2(\mb{R}^n,e^{-\vp})$. Let $u=u_0$. So there exists $u\in L^2(\mb{R}^n,e^{-\vp})$ such that
$\Delta u+au=f$ with $\|u\|^2_\vp\leq c$. The proof is complete.
\ep

\bl\label{lemmaH}
Let $\vp\in C^{4}(\mb{R}^n)$. Then
\begin{align*}
\lt\|\lt(\Delta+a\rt)^*_\vp\phi\rt\|^2_\vp
=\lt\|\lt(\Delta+a\rt)\phi\rt\|^2_\vp+\lt\langle \phi,\Delta\lt(\Delta_\vp^{*}\phi\rt)-\Delta^{*}_\vp
\lt(\Delta\phi\rt)\rt\rangle_\vp, \ \ \forall\phi\in C_0^\infty(\mb{R}^n).
\end{align*}
\el

\bp
Let $\Delta+a=H$. Then $\lt(\Delta+a\rt)^*_\vp=H^*_\vp$.
For $\forall\phi\in C_0^\infty(\mb{R}^n)$, we have
\begin{align}
\lt\|H^*_\vp\phi\rt\|^2_\vp&=\lt\langle H^*_\vp\phi,H^*_\vp\phi\rt\rangle_\vp\nonumber\\
&=\lt\langle \phi,HH^*_\vp\phi\rt\rangle_\vp\nonumber\\
&=\lt\langle \phi,H^*_\vp H\phi\rt\rangle_\vp+\lt\langle \phi,HH^*_\vp\phi-H^*_\vp H\phi\rt\rangle_\vp\nonumber\\
&=\lt\langle H\phi,H\phi\rt\rangle_\vp+\lt\langle \phi,HH^*_\vp\phi-H^*_\vp H\phi\rt\rangle_\vp\nonumber\\
&=\lt\|H\phi\rt\|^2_\vp+\lt\langle \phi,HH^*_\vp\phi-H^*_\vp H\phi\rt\rangle_\vp\label{31}
\end{align}
Note that
\begin{align*}
HH^*_\vp\phi&=\lt(\Delta+a\rt)\lt(\Delta+a\rt)^*_\vp\phi\\
&=\lt(\Delta+a\rt)\lt(\Delta^{*}_\vp\phi+a\phi\rt)\\
&=\Delta\lt(\Delta_\vp^{*}\phi\rt)+a\Delta\phi+
a\Delta_\vp^{*}\phi+a^2\phi
\end{align*}
and
\begin{align*}
H^*_\vp H\phi&=\lt(\Delta+a\rt)^*_\vp\lt(\Delta+a\rt)\phi\\
&=\lt(\Delta^{*}_\vp+a\rt)\lt(\Delta\phi+a\phi\rt)\\
&=\Delta^{*}_\vp\lt(\Delta\phi\rt)+a\Delta_\vp^{*}\phi+a\Delta\phi+a^2\phi.
\end{align*}
Then
\begin{equation}\label{32}
HH^*_\vp\phi-H^*_\vp H\phi=\Delta\lt(\Delta_\vp^{*}\phi\rt)-\Delta^{*}_\vp\lt(\Delta\phi\rt).
\end{equation}
So by (\ref{31}) and (\ref{32}), we have
\begin{align*}
\lt\|H^*_\vp\phi\rt\|^2_\vp=\lt\|H\phi\rt\|^2_\vp+\lt\langle \phi,\Delta\lt(\Delta_\vp^{*}\phi\rt)-\Delta^{*}_\vp\lt(\Delta\phi\rt)\rt\rangle_\vp.
\end{align*}
This lemma is proved.
\ep

\bl\label{lemmak} Let $\vp=|x|^2$. Then for $\forall\phi\in C_0^\infty(\mb{R}^n)$,
\begin{align*}
\lt\ll\phi, \Delta\lt(\Delta_\vp^{*}\phi\rt)-\Delta^{*}_\vp\lt(\Delta\phi\rt)
\rt\rl_\vp=8n\|\phi\|^2_\vp+8\|\nabla\phi\|^2_\vp.
\end{align*}
\el

\bp
For $\forall\phi\in C_0^\infty(\mb{R}^n)$, by the definition of $\Delta_\vp^{*}$ and the following formula
$$\Delta(\alpha\beta)=\beta\Delta\alpha+\alpha\Delta\beta
+2\nabla\alpha\cdot\nabla\beta,\ \ \forall\alpha,\beta\in C^2(\mb{R}^n),$$ we have
\begin{align}\label{1}
\Delta_\vp^{*}\phi=e^{\vp}\Delta\lt(\phi e^{-\vp}\rt)=\Delta\phi+\phi|\nabla\vp|^2-\phi\Delta\vp-2\nabla\phi\cdot
\nabla\vp.
\end{align}
From (\ref{1})
we have
\begin{align*}
\Delta\lt(\Delta_\vp^{*}\phi\rt)&=\Delta^2\phi+\Delta(\phi|\nabla\vp|^2)
-\Delta(\phi\Delta\vp)-2\Delta(\nabla\phi\cdot
\nabla\vp)\\
&=\Delta^2\phi+\Delta\phi|\nabla\vp|^2+\phi\Delta(|\nabla\vp|^2)+
2\nabla\phi\cdot\nabla(|\nabla\vp|^2)\\
&\ \ \ -\Delta\phi\Delta\vp
-\phi\Delta^2\vp-2\nabla\phi\cdot\nabla(\Delta\vp)
-2\Delta(\nabla\phi\cdot
\nabla\vp)
\end{align*}
and
$$\Delta^*_\vp(\Delta\phi)=\Delta^2\phi+\Delta\phi|\nabla\vp|^2
-\Delta\phi\Delta\vp-2\nabla(\Delta\phi)\cdot\nabla\vp.$$
Then
\begin{align}
\Delta\lt(\Delta_\vp^{*}\phi\rt)-\Delta^{*}_\vp\lt(\Delta\phi\rt)&=
\phi\Delta(|\nabla\vp|^2)+
2\nabla\phi\cdot\nabla(|\nabla\vp|^2)
-\phi\Delta^2\vp\nonumber\\&\ \ \ -2\nabla\phi\cdot\nabla(\Delta\vp)
-2\Delta(\nabla\phi\cdot
\nabla\vp)+2\nabla(\Delta\phi)\cdot\nabla\vp\label{41}.
\end{align}
Let $\vp=|x|^2$. We have $\nabla\vp=2x$, $\Delta\vp=2n$,
$|\nabla\vp|^2=4|x|^2$, $\nabla(|\nabla\vp|^2)=8x$, $\Delta(|\nabla\vp|^2)=8n$.
Then by (\ref{41}) and the following formula
$$\Delta(\nabla\phi\cdot x)=\nabla(\Delta\phi)\cdot x+2\Delta\phi,\ \ \forall\phi\in C_0^\infty(\mb{R}^n),$$ we get
\begin{align*}
\Delta\lt(\Delta_\vp^{*}\phi\rt)-\Delta^{*}_\vp\lt(\Delta\phi\rt)
=8n\phi+16(\nabla\phi\cdot x)-8\Delta\phi.
\end{align*}
Consequently,
\begin{align*}
\lt\ll\phi, \Delta\lt(\Delta_\vp^{*}\phi\rt)-\Delta^{*}_\vp\lt(\Delta\phi\rt)
\rt\rl_\vp&=\lt\ll\phi, 8n\phi+16(\nabla\phi\cdot x)-8\Delta\phi
\rt\rl_\vp\\
&=8n\|\phi\|^2_\vp+8\lt\ll\phi, 2(\nabla\phi\cdot x)-\Delta\phi
\rt\rl_\vp.
\end{align*}
Note, as the key step of the proof, that
\begin{align*}
\lt\ll\phi, 2(\nabla\phi\cdot x)-\Delta\phi
\rt\rl_\vp&=\int_{\mb{R}^n}\phi(2(\nabla\phi\cdot x)-\Delta\phi)e^{-\vp}dx\\
&=\int_{\mb{R}^n}\phi\sum^n_{j=1}
\lt(2x_j\frac{\partial\phi}{\partial x_j}-\frac{\partial^2\phi}{\partial x_j\partial x_j}\rt)e^{-|x|^2}dx\\
&=-\int_{\mb{R}^n}\phi\sum^n_{j=1}\frac{\partial}{\partial x_j}\lt(\frac{\partial\phi}{\partial x_j}e^{-|x|^2}\rt)dx\\
&=-\sum^n_{j=1}\int_{\mb{R}^n}\phi\frac{\partial}{\partial x_j}\lt(\frac{\partial\phi}{\partial x_j}e^{-|x|^2}\rt)dx\\
&=\sum^n_{j=1}\int_{\mb{R}^n}\frac{\partial\phi}{\partial x_j}\lt(\frac{\partial\phi}{\partial x_j}e^{-|x|^2}\rt)dx\\
&=\sum^n_{j=1}\int_{\mb{R}^n}\lt(\frac{\partial\phi}{\partial x_j}\rt)^2e^{-|x|^2}dx\\
&=\int_{\mb{R}^n}|\nabla\phi|^2e^{-|x|^2}dx\\
&=\|\nabla\phi\|^2_\vp.
\end{align*}
Then
\begin{align*}
\lt\ll\phi, \Delta\lt(\Delta_\vp^{*}\phi\rt)-\Delta^{*}_\vp\lt(\Delta\phi\rt)
\rt\rl_\vp=8n\|\phi\|^2_\vp+8\|\nabla\phi\|^2_\vp.
\end{align*}
The lemma is proved.
\ep

\section{Proof of theorems}

The proof of Theorems \ref{thk2}.

\bp
Let $\vp=|x|^2$. By Lemma \ref{lemmaH} and Lemma \ref{lemmak}, we have for $\forall\phi\in C_0^\infty(\mb{R}^n)$,
\begin{align}
\lt\|\lt(\Delta+a\rt)^*_\vp\phi\rt\|^2_\vp
&=\lt\|\lt(\Delta+a\rt)\phi\rt\|^2_\vp+\lt\langle \phi,\Delta\lt(\Delta_\vp^{*}\phi\rt)-\Delta^{*}_\vp
\lt(\Delta\phi\rt)\rt\rangle_\vp\nonumber\\
&\geq\lt\langle \phi,\Delta\lt(\Delta_\vp^{*}\phi\rt)-\Delta^{*}_\vp
\lt(\Delta\phi\rt)\rt\rangle_\vp\nonumber\\
&\geq8n\|\phi\|^2_\vp.\label{34}
\end{align}
By Cauchy-Schwarz inequality and (\ref{34}), we have for $\forall\phi\in C_0^\infty(\mb{R}^n)$,
\begin{align*}
|\langle f,\phi\rangle_\vp|^2
&\leq\lt\|f\rt\|^2_\vp
\lt\|\phi\rt\|^2_\vp\\
&=\lt(\frac{1}{8n}\lt\|f\rt\|^2_\vp\rt)
\lt(8n\lt\|\phi\rt\|^2_\vp\rt)\\
&\leq\lt(\frac{1}{8n}\lt\|f\rt\|^2_\vp\rt)
\lt\|\lt(\Delta+a\rt)^*_\vp\phi\rt\|^2_\vp.
\end{align*}
Let $c=\frac{1}{8n}\lt\|f\rt\|^2_\vp$. Then
$$|\langle f,\phi\rangle_\vp|^2\leq c\lt\|\lt(\Delta+a\rt)^*_\vp\phi\rt\|^2_\vp,\ \ \forall\phi\in C_0^\infty(\mb{R}^n).$$
By Lemma \ref{lemmaifif}, there exists a global weak solution $u\in L^2(\mb{R}^n,e^{-\vp})$ solving the equation
$$\Delta u+au=f$$ in $\mb{R}^n$
with the norm estimate
$$\|u\|^2_\vp\leq c,$$
i.e.,
$$\Delta u+au=f \ \ {\it with} \ \ \int_{\mb{R}^n} u^2e^{-|x|^2}dx\leq
\frac{1}{8n}\int_{\mb{R}^n} f^2e^{-|x|^2}dx.$$
The proof is complete.
\ep

\bt\label{thwj}
There exists a bounded operator $Q:L^2(\mb{R}^n,e^{-|x|^2})\rightarrow L^2(\mb{R}^n,e^{-|x|^2})$ such that
$$\lt(\Delta+a\rt)Q=I \ \ {\it with} \ \ \|Q\|\leq\frac{1}{\sqrt{8n}},$$
where $\|Q\|$ is the norm of $Q$ in $L^2(\mb{R}^n,e^{-|x|^2})$.
\et

\bp
Let $\vp=|x|^2$. For each $f\in L^2(\mb{R}^n,e^{-\vp})$, from Theorem \ref{thk2},
there exists $u\in L^2(\mb{R}^n,e^{-\vp})$ such that
$$\lt(\Delta+a\rt)u=f \ \ {\it with} \ \ \|u\|_\vp\leq\frac{1}{\sqrt{8n}}\lt\|f\rt\|_\vp.$$
Denote this $u$ by $Q(f)$. Then $Q(f)$ satisfies $$\lt(\Delta+a\rt)Q(f)=f \ \ {\it with} \ \ \|Q(f)\|_\vp\leq\frac{1}{\sqrt{8n}}\lt\|f\rt\|_\vp.$$
Note that $f$ is arbitrary in $L^2(\mb{R}^n,e^{-\vp})$. So $Q:L^2(\mb{R}^n,e^{-\vp})\rightarrow L^2(\mb{R}^n,e^{-\vp})$ is a bounded operator such that
$$\lt(\Delta+a\rt)Q=I \ \ {\it with} \ \ \|Q\|\leq\frac{1}{\sqrt{8n}}.$$ The proof is complete.
\ep

\section{Further remarks}

\noindent\textbf{Remark 1.}
Given $\lambda>0$ and $x_0\in\mb{R}^n$, for the weight $\vp=\lambda|x-x_0|^2$, we obtain the following corollary from Theorem \ref{thk2}.

\bc\label{thzy}
For each $f\in L^2(\mb{R}^n,e^{-\lambda|x-x_0|^2})$, there exists a weak solution
$u\in L^2(\mb{R}^n,e^{-\lambda|x-x_0|^2})$ solving the equation
$$\Delta u+au=f$$
with the norm estimate $$\int_{\mb{R}^n} u^2e^{-\lambda|x-x_0|^2}dx\leq
\frac{1}{8n\lambda^2}\int_{\mb{R}^n} f^2e^{-\lambda|x-x_0|^2}dx.$$
\ec

\bp
From $f\in L^2(\mb{R}^n,e^{-\lambda|x-x_0|^2})$, we have
\begin{equation}\label{101}
\int_{\mb{R}^n} f^2(x)e^{-\lambda|x-x_0|^2}dx<+\infty.
\end{equation}
Let $x=\frac{y}{\sqrt{\lambda}}+x_0$ and $g(y)=f(x)= f\lt(\frac{y}{\sqrt{\lambda}}+x_0\rt)$. Then by (\ref{101}), we have
\begin{equation*}
\frac{1}{\lt(\sqrt{\lambda}\rt)^n}\int_{\mb{R}^n} g^2(y)e^{-|y|^2}dy<+\infty,
\end{equation*}
which implies that $g\in L^2(\mb{R}^n,e^{-|y|^2})$. For $g$, applying Theorem \ref{thk2} with $a$ replaced by $\frac{a}{\lambda}$, there exists a weak solution $v\in L^2(\mb{R}^n,e^{-|y|^2})$ solving the equation
\begin{equation}\label{103}
\Delta v(y)+\frac{a}{\lambda}v(y)=g(y)
\end{equation}
in $\mb{R}^n$ with the norm estimate
\begin{equation}\label{104}
\int_{\mb{R}^n} v^2(y)e^{-|y|^2}dy\leq
\frac{1}{8n}\int_{\mb{R}^n} g^2(y)e^{-|y|^2}dy.
\end{equation}
Note that $y=\sqrt{\lambda}(x-x_0)$ and $g(y)=f(x)$. Let $u(x)=\frac{1}{\lambda}v(y)=\frac{1}{\lambda}v\lt(\sqrt{\lambda}(x-x_0)\rt)$. Then (\ref{103}) and (\ref{104}) can be rewritten by
\begin{equation}\label{109}
\Delta u(x) +au(x)
=f(x)
\end{equation}
\begin{equation}\label{110}
\int_{\mb{R}^n} u^2(x)e^{-\lambda|x-x_0|^2}dx\leq
\frac{1}{8n\lambda^2}\int_{\mb{R}^n} {f^2(x)}e^{-\lambda|x-x_0|^2}dx.
\end{equation}
(\ref{110}) implies that $u\in L^2(\mb{R}^n,e^{-\lambda|x-x_0|^2})$. Then by (\ref{109}) and (\ref{110}), the proof is complete.
\ep

\noindent\textbf{Remark 2.}
From Corollary \ref{thzy}, we can obtain the following corollary, which shows that for any choice of $a$, the differential operator $\Delta+a$ has a bounded right inverse in $L^2(U)$, provided $U$ is a bounded open set. This result should be well-known, and we could not locate a reference (for a very close related result, see Theorem 6 in page 324 in \cite{2}).
\bc\label{thzwj1}
Let $U\subset\mb{R}^n$ be any bounded open set. For each $f\in L^2(U)$, there exists a weak solution $u\in L^2(U)$ solving the equation
$$\Delta u+au=f$$
with the norm estimate $\|u\|_{L^2(U)}\leq c\|f\|_{L^2(U)}$,
where the constant $c$ depends only on the diameter of $U$.
\ec

\bp
Let $x_0\in U$. For given $f\in L^2(U)$, extending $f$ to zero on $\mb{R}^n\setminus U$, we have
\begin{align*}
\tilde{f}=\left\{
\begin{array}{ccc}
f,       &      & {x\in U}\\
0,     &      & {x\in \mb{R}^n\setminus U.}
\end{array} \right.
\end{align*}
Then $\tilde{f}\in L^2(\mb{R}^n)\subset L^2(\mb{R}^n,e^{-|x-x_0|^2})$. From Corollary \ref{thzy},
there exists $\tilde{u}\in L^2(\mb{R}^n,e^{-|x-x_0|^2})$ such that
$$\Delta\tilde{u}+a\tilde{u}=\tilde{f} \ \ {\it with} \ \ \int_{\mb{R}^n}\tilde{u}^2e^{-|x-x_0|^2}dx
\leq\frac{1}{8n}\int_{\mb{R}^n}\tilde{f}^2e^{-|x-x_0|^2}dx.$$ Then
\begin{align*}
\int_{\mb{R}^n}\tilde{u}^2e^{-|x-x_0|^2}dx
\leq\frac{1}{8n}\int_{\mb{R}^n}\tilde{f}^2dx
=\frac{1}{8n}\int_Uf^2dx.
\end{align*}
Note that
\begin{align*}
\int_{\mb{R}^n}\tilde{u}^2e^{-|x-x_0|^2}dx
\geq\int_U\tilde{u}^2e^{-|x-x_0|^2}dx
\geq\int_U\tilde{u}^2e^{-|U|^2}dx
=e^{-|U|^2}\int_U\tilde{u}^2dx,
\end{align*}
where $|U|$ is the diameter of $U$.
Therefore,
\begin{align*}
e^{-|U|^2}\int_U\tilde{u}^2dx\leq \frac{1}{8n}\int_Uf^2dx, \ \ i.e.,\ \ \int_U\tilde{u}^2dx\leq\frac{e^{|U|^2}}{8n}\int_Uf^2dx.
\end{align*}
Restricting $\tilde{u}$ on $U$ to get $u$,
then
$$\Delta u+au=f \ \ {\it with} \ \ \int_Uu^2dx\leq  \frac{e^{|U|^2}}{8n}\int_Uf^2dx.$$
Note that $u\in L^2(U)$ and let $c=\sqrt{\frac{e^{|U|^2}}{8n}}$. Then the proof is complete.
\ep

\noindent\textbf{Remark 3.} When $f\in L^2(\mb{R}^n)$, the solutins of $\Delta u=f$ are not necessary in $L^2(\mb{R}^n)$. For example: $n=1$,
\begin{align*}
f(x)=\left\{
\begin{array}{ccc}
\frac{1}{x},       &      & {x\geq1,}\\
x,     &      & {0<x<1,}\\
0,     &      & {x\leq0,}
\end{array} \right.
\end{align*}
$$u(x)=\int^x_0(x-t)f(t)dt+c_1x+c_2
=-\frac{x}{2}+xlnx+\frac{2}{3}+c_1x+c_2,\ \ x\geq1,$$
where $c_1$ and $c_2$ are arbitrary real constants. It is easy to see $u\not\in L^2(\mb{R})$.


\begin{thebibliography}{XX}

\bibitem{2} L. C. Evans, Partial differential equations, Second edition (2010).

\bibitem{1} L. H\"{o}rmander, $L^2$ estimates and existence theorems for the $\overline{\partial}$ operator, Acta Mathematica, 113 (1965), 89-152.



\end{thebibliography}
\end{document}